\documentclass[10pt,oneside]{jm}
\usepackage{amsfonts,amssymb,amsmath,graphicx}

\newcommand{\abs}{\vskip 0.5em\noindent}

\newcommand{\AbsT}[1]{\vskip 0.5em\noindent {\bf #1} }
\newcommand{\AbsC}[1]{\begin{center}\bf #1 \end{center} \vskip -0.5em}

\newcommand{\Pf}{\AbsT{Proof.}}

\newcommand{\Rem}{\AbsT{Remark.}}

\newcommand{\Prop}{\AbsT{Proposition.}}

\newcommand{\Thm}{\AbsT{Theorem.}}

\newcommand{\Cor}{\AbsT{Corollary.}}

\newcommand{\bfi}{\noindent {\bf i)} }
\newcommand{\bfii}{\noindent {\bf ii)} }
\newcommand{\bfiii}{\noindent {\bf iii)} }

\newcommand{\bfa}{\noindent {\bf a)} }
\newcommand{\bfb}{\noindent {\bf b)} }

\newcommand{\ie}{i.~e.\hskip 0.5em}

\newcommand{\QED}{\hfill $\sharp$}

\newcommand{\N}{\mathbb N}

\newcommand{\R}{\mathbb R}
\newcommand{\C}{\mathbb C}

\newcommand{\cA}{\mathcal A}

\newcommand{\cC}{\mathcal C}
\newcommand{\cD}{\mathcal D}

\newcommand{\cF}{\mathcal F}

\newcommand{\cP}{\mathcal P}

\newcommand{\cQ}{\mathcal Q}

\newcommand{\cS}{\mathcal S}

\newcommand{\cX}{\mathcal X}

\newcommand{\al}{\alpha}

\newcommand{\bt}{\beta}

\newcommand{\dt}{\delta}

\newcommand{\lb}{\lambda}

\newcommand{\Aut}{\text{Aut}}

\newcommand{\diag}{\text{diag}}

\newcommand{\GL}{\text{GL}}

\newcommand{\Hom}{\text{Hom}}

\newcommand{\id}{\text{id}}

\newcommand{\tr}{\text{tr}}

\newcommand{\baraut}{\ov{\text{\hspace{0.5em}\rule{0em}{0.55em}}}}

\newcommand{\ld}{,\ldots\hskip0em ,}
\newcommand{\lr}[1]{\langle #1\rangle}
\newcommand{\ov}[1]{\overline{#1}}

\newcommand{\fl}[1]{\lfloor #1\rfloor}

\newcommand{\wti}[1]{\widetilde{#1}}

\newcommand{\mt}{\mapsto}

\newcommand{\ra}{\rightarrow}

\newcommand{\cn}{\colon}
\newcommand{\dcup}{\stackrel{.}{\cup}}

\newcommand{\sseq}{\subseteq}
\newcommand{\smin}{\setminus}

\newcommand{\ity}{\infty}

\newcommand{\bsl}{\backslash}
\newcommand{\tm}{\times}

\begin{document}
\raggedbottom
\pagestyle{myheadings}
\markboth{}{}
\thispagestyle{empty}
\setcounter{page}{1}
\begin{center} \Large\bf
Symmetric generation of Coxeter groups
\vspace*{1em} \\
\large\rm
Ben Fairbairn and J\"urgen M\"uller
\vspace*{1em} \\
\end{center}

\begin{abstract} \noindent
We provide involutory symmetric generating sets of finitely generated
Coxeter groups, fulfilling a suitable finiteness condition,
which in particular is fulfilled in the finite, affine and 
compact hyperbolic cases. 
\end{abstract}

\AbsC{Introduction}

\abs
Presentations of groups having certain types of symmetry properties
are being considered at least since Coxeter's work \cite{Cox}
in 1959. The notion of symmetric generation we are interested in here
has been introduced by Curtis \cite{C2} in 1992.
Since then it has gained continuing interest, in particular in the theory 
of sporadic finite simple groups, see Curtis's recent book \cite{C3} 
and the references in there, typically leading to `good' presentations
providing a practical means of computing efficiently 
in the groups under consideration \cite{CF}.

\abs
We recall the definition of {\bf symmetric generation} 
as we use it in this note:
Let $k^{\ast m}$ be the free product of
$m\in\N$ copies of the cyclic groups $C_k$ of order $k\geq 2$, where 
as usual the latter are often abbreviated by just writing $k$.
Then $\Aut(C_k)$ acts on any of the free factors,
and together with the full symmetric group $\cS_m$ permuting 
the free factors, this generates the group $\Aut(C_k)\wr\cS_m$
of monomial automorphisms of $k^{\ast m}$; for $k=2$ this
of course simplifies to purely permutational automorphisms.
Given a {\bf control group} $\cC\leq\Aut(C_k)\wr\cS_m$ of 
monomial automorphisms acting transitively on the free factors, the 
semidirect product $k^{\ast m}\cn\cC$ is called a {\bf progenitor}.

\abs
An epimorphic image $G$ of a progenitor is called
{\bf (strictly) symmetrically generated} if 
\bfi the image of $\cC$ in $G$ is isomorphic to $\cC$,
\bfii the free factors map to $m$ cyclic subgroups of $G$ 
of order $k$ having mutually trivial intersections, and
\bfiii the latter generate $G$.
If only the first two conditions are fulfilled, but the third 
is not, then $G$ is called {\bf weakly} symmetrically generated;
this less restrictive notion indeed
leads to certain interesting examples \cite{B2,B1}.

\abs
The aim of this note is to provide involutory symmetric 
generating sets of finitely generated, but not necessarily finite
Coxeter groups. Given such a group,
for any of its maximal parabolic subgroups fulfilling
a suitable finiteness condition we obtain an associated 
weakly symmetric generating set, and we settle the question 
when it is strict. 
In particular, as far as irreducible Coxeter groups
are concerned, the necessary conditions are fulfilled 
for all maximal parabolic subgroups of finite or affine Coxeter groups,
\ie the tame cases in the sense of \cite[Sect.1.11]{L},
and of compact hyperbolic Coxeter groups \cite[Sect.6.8]{H},
while for most non-compact hyperbolic Coxeter groups \cite[Sect.6.9]{H}
there is at least one suitable maximal parabolic subgroup.

\abs
This note is motivated by \cite{F}, and gives a general 
explanation of the observations made there. As it turns out, 
Coxeter groups fit extremely nicely into the concept of 
symmetric generation. This to our knowledge has
escaped notice so far, despite the fact that
the symmetric group being generated by the
adjacent transpositions is the archetypical example
of symmetric generation \cite[Exc.3.3(1)]{C3}.
Moreover, not merely finite Coxeter groups but just finitely generated 
Coxeter groups are the natural appropriate class  
to state our results for, and hence we try to minimise our finiteness
assumptions as far as possible.
To our knowledge this is the first class of possibly
infinite groups for which symmetric generation in the above sense is
being considered at all. It remains to be seen whether the
observations made in this note will be helpful to handle Coxeter groups 
in practice and to gain further insights into them.

\abs
This note is organised as follows: 
In Section 1 we introduce the progenitors,
in Section 2 we prove the main result on symmetric generation,
in Section 3 we examine the permutation action of the control group,
and in Section 4 we consider various explicit examples.
We assume the reader to be familiar with presentations of groups,
a basic reference being \cite{J}, and with Coxeter groups and 
related concepts, as exposed in \cite{B,H}. 
Our notation for groups is as used in \cite{Atlas}.
We consider right actions throughout.

\AbsC{1. The progenitor}

\abs
Let $W$ be a finitely generated, not necessarily finite
Coxeter group with distinguished 
generators $\{s_1\ld s_n\}\sseq W$, for some $n\in\N$,
and associated Coxeter integers $m_{ij}\in\N\dcup\{\ity\}$ such that
$m_{ii}=1$ and $m_{ij}=m_{ji}>1$, 
for $i,j\in\{1\ld n\}$ such that $i\neq j$.
This gives rise to the Coxeter presentation
$$ W\cong\cF(S_1\ld S_n)/\lr{\!\lr{(S_iS_j)^{m_{ij}};
             i,j\in\{1\ld n\}}\!} ,$$
where $\cF(S_1\ld S_n)$ denotes the free group with generators
$S_1\ld S_n$, where $\lr{\!\lr{\cdot}\!}$ denotes normal closure,
and where $(S_iS_j)^{\ity}=1$ is interpreted as being no relation at all.
We have $W=W_1\tm\cdots\tm W_t$,
where the direct factors are the parabolic subgroups associated with
the connected components of the Dynkin graph associated with $W$; 
we may assume that $s_n\in W_t$.
Let $W':=\lr{s_1\ld s_{n-1}}<W$, which is a maximal parabolic subgroup
and hence by \cite[Prop.9.5]{L} has the Coxeter presentation 
$$ W'\cong\cF(S_1\ld S_{n-1})/\lr{\!\lr{(S_iS_j)^{m_{ij}};
              i,j\in\{1\ld n-1\}}\!} ,$$
obtained from the presentation of $W$ 
by leaving out the relations involving $S_n$.
Letting $W'_t:=W_t\cap W'$ we have $W'=W_1\tm\cdots\tm W_{t-1}\tm W'_t$.

\abs
Let $(s_n)^{W'}\sseq W$ be the set of $W'$-conjugates of $s_n$,
and let $C_{W'}(s_n)\leq W'$ be the centraliser of $s_n$ in $W'$.
{\bf We assume that $C_{W'}(s_n)\leq W'$ has finite index.}
Since $W_1\tm\cdots\tm W_{t-1}\leq C_{W'}(s_n)$
this in particular is fulfilled whenever $W'_t$ is finite.
Let $m:=[W'\cn C_{W'}(s_n)]=|(s_n)^{W'}|\in\N$, 
and let $\{w_1\ld w_m\}\sseq W'$ 
be a set of representatives of the right cosets $C_{W'}(s_n)\bsl W'$ 
of $C_{W'}(s_n)$ in $W'$, where $w_1:=1$, yielding the bijection 
$C_{W'}(s_n)\bsl W'\ra (s_n)^{W'}\cn w_k\mt (s_n)^{w_k}=:t_k$
for $k\in\{1\ld m\}$, where $t_1=s_n$.
Hence the action of $W'$ on $C_{W'}(s_n)\bsl W'$
translates to the transitive permutation representation $\pi\cn W'\ra\cS_m$,
where $\pi(w)\in\cS_m$ is given by 
$t_{k^{\pi(w)}}=(t_k)^w$ for $k\in\{1\ld m\}$.

\abs
We consider the finitely presented group
$$ \begin{array}{rl}
\cP:=\cF(S_1\ld S_{n-1},T_1\ld T_m)/\lr{\!\lr{\hspace*{-0.8em}&
(S_iS_j)^{m_{ij}},\,\,\, T_1^2,\,\,\, (T_k)^{S_i}T_{k^{\pi(s_i)}}^{-1}; \\
&i,j\in\{1\ld n-1\},k\in\{1\ld m\}}\!}. \\
\end{array} $$
Since $\pi$ is transitive, $\cP$ is generated by the images of
$S_1\ld S_{n-1},T_1$, and the relations $T_k^2=1$ for $k\in\{2\ld m\}$
follow from the relations given.
This implies that $\cP\cong 2^{\ast m}\cn W'$,
where $W'$ acts by permuting the 
free factors of $2^{\ast m}$ via $\pi$.
Thus $\cP$ is a progenitor with control group $W'$.

\AbsC{2. Symmetric generation}

\Thm
We still assume that $m:=[W'\cn C_{W'}(s_n)]$ is finite. 
Then $W$ has the finite presentation 
$$ \cQ:=\cP/\lr{\!\lr{(S_iT_1)^{m_{in}};
                      i\in\{1\ld n-1\}\text{ such that }m_{in}>2}\!} .$$

\Pf
Let $\cF:=\cF(S_1\ld S_{n-1},T_1\ld T_m)$,
and let $\baraut\cn\cF\ra\cQ$ be the natural epimorphism.
Letting $S_i\mt s_i$ for $i\in\{1\ld n-1\}$, and 
$T_k\mt (s_n)^{w_k}$ for $k\in\{1\ld m\}$
defines a homomorphism $\cF\ra W$. 
In particular we have $T_1\mt s_n$, hence this is an epimorphism.
Since the defining relations of $\cP$,
and the relations $(S_iT_1)^{m_{in}}=1$ for $i\in\{1\ld n-1\}$
are fulfilled by the relevant images in $W$, 
this induces epimorphisms $\cP\ra W$ and 
$\al\cn\cQ\ra W\cn\ov{S}_i\mt s_i,\ov{T}_k\mt (s_n)^{w_k}=t_k$
for $i\in\{1\ld n-1\}$ and $k\in\{1\ld m\}$.

\abs
As for the converse, we have an epimorphism $\cF(S_1\ld S_n)\ra\cQ$
by letting $S_i\mt\ov{S}_i$ for $i\in\{1\ld n-1\}$, and $S_n\mt\ov{T}_1$ .
Since the defining relations of $W$ involving $S_iS_j$ for 
$i,j\in\{1\ld n-1\}$, and involving $S_iS_n$ for $i\in\{1\ld n-1\}$ 
such that $m_{in}>2$ are by construction fulfilled by the relevant 
images in $\cQ$, only those involving $S_iS_n$ where $m_{in}=2$ 
remain to be dealt with:
For $i\in\{1\ld n-1\}$ such that $m_{in}=2$ we have
$s_i\in C_{W'}(s_n)$, hence
$t_{1^{\pi(s_i)}}=t_1^{s_i}=(s_n)^{s_i}=s_n=t_1$, 
and thus $1^{\pi(s_i)}=1$,
implying that the relation $(T_1)^{S_i}=T_{1^{\pi(s_i)}}=T_1$,
being equivalent to $(S_iT_1)^2=1$, 
is already included in the defining relations of $\cP$.
Hence we have an epimorphism 
$\bt\cn W\ra\cQ\cn s_i\mt\ov{S}_i, s_n\mt\ov{T}_1$
for $i\in\{1\ld n-1\}$. 
Thus $\al$ and $\bt$ are mutually inverse isomorphisms.
\QED

\abs
The above proof also shows that
the image of $W'\leq\cP$ in $\cQ$ is
isomorphic to $W'$, and the images of $T_1\ld T_m$ in $\cQ$
are pairwise distinct involutions,
\ie $\cQ$ is weakly symmetrically generated.
Hence it remains to determine when precisely $\cQ$ is 
strictly symmetrically generated, \ie
when $W$ is generated by $t_1\ld t_m$. 
We proceed as follows:

\abs
We first describe when distinguished generators $s_i$ and $s_j$ 
are conjugate in $W$:
Let $\sim$ be the finest equivalence relation on $\{s_1\ld s_n\}$
such that $s_i\sim s_j$ whenever $m_{ij}$ is odd.
The relation $(s_is_j)^{m_{ij}}=1$ shows that 
$\lr{s_i,s_j}\cong\cD_{2m_{ij}}$ is isomorphic to
the dihedral group of order $2m_{ij}$.
Hence if $m_{ij}$ is odd, then $s_i$ and $s_j$ are already 
conjugate in $\lr{s_i,s_j}$, and thus in $W$.
This shows that if $s_i\sim s_j$, then $s_i$ and $s_j$ 
are conjugate in $W$.
Conversely, we pick a $\sim$-equivalence class $[s]$, and 
let $\lb_{[s]}\cn\{s_1\ld s_n\}\ra\{\pm 1\}\leq\C^\ast$
be defined by $\lb_{[s]}(s_i):=-1$ if and only if $s_i\sim s$.
Then the defining relations of $W$ show that 
$\lb_{[s]}$ extends to a linear character of $W$.
This shows that if $s_i\not\sim s_j$, then $s_i$ and $s_j$
are not conjugate in $W$.

\abs
Let $[W,W]\leq W$ denote the derived subgroup of $W$.
Since $W$ is generated by involutions, $W/[W,W]$ is $2$-elementary abelian.
The above argument on linear characters shows that $W/[W,W]\cong 2^r$,
where $r$ is the number of $\sim$-equivalence classes,
and the direct factors are generated by images of representatives
of the $\sim$-equivalence classes.

\Prop
$W$ is strictly symmetrically generated, 
\ie $W=\lr{t_1\ld t_m}$, if and only if $r=1$,
\ie if and only if the subgraph of the Dynkin graph of $W$
consisting of the edges labelled by odd $m_{ij}$ is connected.

\Pf
Since $t_1\ld t_m$ are all conjugate in $W$ to $s_n$, we conclude
that $W=\lr{t_1\ld t_m}$ implies $r=1$. If conversely $r=1$, then
after possibly reordering $s_1\ld s_{n-1}$ we may assume that for any 
$i\in\{1\ld n-1\}$ there is $j\in\{i+1\ld n\}$ such that $m_{ij}$ is odd. 
We for $i\in\{n,n-1\ld 1\}$ successively show that $s_i\in\lr{t_1\ld t_m}$:
We have $s_n=t_1$ anyway. For $i<n$ let $j$ be as above.
Hence there are $l\in\N$ and $j_k\in\{1\ld m\}$
such that $s_j=t_{j_1}\cdot t_{j_2}\cdot\cdots\cdot t_{j_l}$.
We have 
$s_j^{s_i}
=t_{j_1}^{s_i}\cdot\cdots\cdot t_{j_l}^{s_i}
=t_{j_1^{\pi(s_i)}}\cdot\cdots\cdot t_{j_l^{\pi(s_i)}}
\in\lr{t_1\ld t_m}$.
Since $m_{ij}$ is odd, we conclude that $s_j\neq s_j^{s_i}$ 
are distinct involutions in $\lr{s_i,s_j}\cong\cD_{2m_{ij}}$, 
hence $\lr{s_i,s_j}=\lr{s_j,s_j^{s_i}}$ implies
$s_i\in\lr{s_j,s_j^{s_i}}\leq\lr{t_1\ld t_m}$.
\QED

\AbsC{3. The representation $\pi$}

\abs
To describe the permutation representation $\pi\cn W'\ra\cS_m$ 
we determine $C_{W'}(s_n)$. 
{\bf In this section no finiteness assumption is needed.}

\Thm
Let $\cX:=\{i\in\{1\ld n-1\};m_{in}=2\}$, and let $W_\cX\leq W'\leq W$
be the associated parabolic subgroup. Then we have $C_{W'}(s_n)=W_\cX$.

\Pf
For $i\in\cX$ we have $s_i\in C_{W'}(s_n)$, and thus 
$W_\cX\leq C_{W'}(s_n)$. To show the converse,
following \cite[Sect.V.4]{B} 
let $E$ be an $\R$-vector space having $\R$-basis $\Pi:=\{\al_1\ld\al_n\}$,
and carrying the symmetric $\R$-bilinear form $B\cn E\tm E\ra\R$ 
given by $b_{ij}:=B(\al_i,\al_j):=-\cos(\frac{\pi}{m_{ij}})$ 
for $i,j\in\{1\ld n\}$; we have $b_{ii}=1$, and for $i\neq j$ we have
$b_{ij}=0$ if and only if $m_{ij}=2$, and $b_{ij}<0$ otherwise.
The geometric action $\rho$ of $W$ on $E$ is given by the 
following representing matrices with respect to $\Pi$,
where we omit zero entries:
$$ \rho(s_i)=\left[ \begin{array}{rrrrrrr}
 1 & & & -2b_{1i} & & & \\
 & \ddots & & \vdots & & & \\
 & & 1 & -2b_{i-1,i} & & & \\
 & & & -1 & & & \\
 & & & -2b_{i+1,i} & 1 & & \\
 & & & \vdots & & \ddots & \\
 & & & -2b_{ni} & & & 1 \\
\end{array} \right]\in\R^{n\tm n} $$
Then $B$ is invariant with respect to $\rho$, 
and by \cite[Cor.V.4.4.2]{B} $\rho$ is faithful. 

\abs
The space $\lr{\al_1\ld\al_{n-1}}_\R\leq E$ is $W'$-invariant,
and $W'$ acts trivially on the quotient space $E/\lr{\al_1\ld\al_{n-1}}_\R$,
hence for any $w\in W'$ we have 
$w\cn\al_n\mt\sum_{i=1}^{n-1}a_i\al_i+\al_n$ for some $a_i\in\R$.
Moreover, the eigenspace of $\rho(s_n)$ with respect to the
eigenvalue $-1$ is given as $E_-(s_n)=\lr{\al_n}_\R$.
Hence any element of $C_W(s_n)$ maps $\lr{\al_n}_\R$ to itself,
and thus from the above we conclude that any element of $C_{W'}(s_n)$
fixes $\al_n$.

\abs
Let $E^\ast:=\Hom_\R(E,\R)$ be the dual space associated with $E$, 
with $\R$-basis $\Pi^\ast:=\{\al^\ast_1\ld\al^\ast_n\}$ dual to $\Pi$.
Then $W$ acts on $E^\ast$ by the contragredient action $\rho^\ast$;
in particular we for $\al\in E$ and $\al^\ast\in E^\ast$ and $w\in W$ have
$\al^w\cdot(\al^\ast)^w=\al\cdot\al^\ast\in\R$,
where we leave out the symbols $\rho$ and $\rho^\ast$
and just write the action of elements as superscripts.
Since the $s_i\in W$ have order $2$, 
we get the following representing matrices
$\rho^\ast(s_i)=\rho(s_i)^{\tr}\in\R^{n\tm n}$
with respect to $\Pi^\ast$: 
$$ \rho^\ast(s_i)=\left[ \begin{array}{rrrrrrr}
 1 & & & & & & \\
 & \ddots & & & & & \\
 & & 1 & & & & \\
 -2b_{i1} & \cdots & -2b_{i,i-1} & -1 & -2b_{i,i+1} & 
            \cdots & -2b_{in} \\
 & & & & 1 & & \\
 & & & & & \ddots & \\
 & & & & & & 1 \\
\end{array} \right] $$

\abs
The space $\lr{\al^\ast_n}_\R\leq E^\ast$ is $W'$-invariant,
and the $W'$-action on the quotient space $E^\ast/\lr{\al^\ast_n}_\R$ 
coincides with the contragredient geometric action of the 
parabolic subgroup $W'$. Hence we may identify 
$E^\ast/\lr{\al^\ast_n}_\R$ with the contragredient geometric space 
$E^{\prime\ast}$ of $W'$, the canonical $\R$-basis of the latter 
being identified with $\{\ov{\al}^\ast_1\ld\ov{\al}^\ast_{n-1}\}$,
where $\baraut\cn E^\ast\ra E^\ast/\lr{\al^\ast_n}_\R$ is the natural map.

\abs
The eigenspace of $\rho^\ast(s_n)$ with respect to the
eigenvalue $-1$ is given as 
$E^\ast_-(s_n)=\ker(\rho^\ast(s_n)+\id_{E^\ast})=\lr{\bt^\ast}_\R$,
where $\bt^\ast:=\sum_{i=1}^{n-1}b_{in}\al_i^\ast+\al_n^\ast$.
Hence any element of $C_W(s_n)$ maps $\lr{\bt^\ast}_\R$ to itself,
thus for $w\in C_{W'}(s_n)$ there is $a\in\R$
such that $\bt^{\ast w}=a\bt^\ast$. From $\al_n^w=\al_n$ we get
$1=\al_n\cdot\bt^\ast=\al_n^{w^{-1}}\cdot\bt^\ast
=\al_n\cdot\bt^{\ast w}=\al_n\cdot(a\bt^\ast)=a$.
Hence any element of $C_{W'}(s_n)$ fixes $\bt^\ast$, and thus
$\ov{\bt}^\ast=\sum_{i=1}^{n-1}b_{in}\ov{\al}_i^\ast\in E^{\prime\ast}$
as well. Considering the cone 
$$ C'_\cX:=\left(\bigcap_{i\in\cX}H'_i\right)\cap
   \left(\bigcap_{i\in\{1\ld n-1\}\smin\cX}A'_i\right)\sseq E^{\prime\ast},$$
where
$H'_i:=\{\al^\ast\in E^{\prime\ast};\al_i\cdot\al^\ast=0\}$ and
$A'_i:=\{\al^\ast\in E^{\prime\ast};\al_i\cdot\al^\ast>0\}$,
we from $\cX=\{i\in\{1\ld n-1\};b_{in}=0\}$ get $-\ov{\bt}^\ast\in C'_\cX$. 
Thus the isotropy group of $\ov{\bt}^\ast$ in $W'$ by \cite[Cor.V.4.6]{B} 
coincides with $W_\cX$, hence we have $C_{W'}(s_n)\leq W_\cX$.
\QED

\Cor
The above proof also shows the following:
The eigenspace of $\rho(s_n)$ with respect to the
eigenvalue $1$ is given as $E_+(s_n)=\lr{\al_n}_\R^\perp$,
where the latter denotes the orthogonal complement
of $\lr{\al_n}_\R$ with respect to $B$.
Since $\rho$ is faithful, $C_W(s_n)$ is 
the set of all elements of $W$ preserving the
vector space decomposition $E=\lr{\al_n}_\R\oplus\lr{\al_n}_\R^\perp$,
and hence we have $C_{W'}(s_n)=\{w\in W';\al_n^w=\al_n\}$.

\abs
Similarly, the eigenspace of $\rho^\ast(s_n)$ with respect to the
eigenvalue $1$ is given as
$E^\ast_+(s_n)=\lr{\al_1^\ast\ld\al_{n-1}^\ast}_\R$.
Hence $C_W(s_n)$ is the set of all elements of $W$
preserving the vector space decomposition
$E^\ast=\lr{\bt^\ast}_\R\oplus\lr{\al_1^\ast\ld\al_{n-1}^\ast}_\R$.

\Rem
We more closely consider the representation $\pi$
in the finite and affine cases.
By factoring out $W_1\tm\cdots\tm W_{t-1}\unlhd C_{W'}(s_n)$
we may assume that $t=1$.
We keep the notation of the above theorem and its proof.

\abs\bfa
Let $W$ be finite. Then a simpler proof to show
$C_{W'}(s_n)\leq W_\cX$ is as follows:
By \cite[Thm.V.4.8.2]{B} $B$ is positive definite. 
Hence we may identify $E$ and $E^\ast$ via $B$,
and $\Pi$ can be considered as a set of fundamental roots in the
root system $\Phi:=\lr{\Pi}_\R\sseq E$. By \cite[Thm.1.12]{H}
$C_{W'}(s_n)=\{w\in W';\al_n^w=\al_n\}$ is generated by the 
reflections it contains. Hence we may assume that $w\in C_{W'}(s_n)$
is a reflection. From \cite[Prop.1.10, Prop.1.14]{H}
we conclude that $w=s_\bt\in W'$ is the reflection
associated with a root $\bt\in\lr{\al_1\ld\al_{n-1}}_\R\cap\Phi$.
We may assume that
$\bt=\sum_{i=1}^{n-1}b_i\al_i$ is positive, \ie we have 
$b_i\geq 0$ for all $i$.
Thus from $B(\bt,\al_n)=0$, and $b_{in}\leq 0$ for all $i\in\{1\ld n-1\}$,
we infer $b_i=0$ whenever $b_{in}<0$, hence 
$\bt\in\lr{\cX}_\R\cap\Phi$ and thus $s_\bt\in W_\cX$.
\QED

\abs\bfb
Let $W$ be affine, hence $W'<W$ is finite and integral.
We consider the particular case where $s_n$ is chosen such that 
the Dynkin graph of $W$ is the completion, 
in the sense of \cite[Sect.VI.4.3]{B}, of the Dynkin graph of $W'$.
We derive another description of $C_{W'}(s_n)=W_\cX$:

\abs
Let $E'$ be the $\R$-vector space underlying the geometric 
representation of $W'$, having a set $\Pi'=\{\al'_1\ld\al'_{n-1}\}$
of fundamental roots in the root system $\Phi':=\lr{\Pi'}_\R\sseq E'$
as an $\R$-basis, 
and carrying the invariant symmetric $\R$-bilinear form $B'$.
Let $-\al'_0\in\Phi'$ be the highest positive root;
if $\Phi'$ contains roots of different lengths, 
then $\al'_0$ is a long root.
We have $B'(\al'_i,\al'_0)=B(\al_i,\al_n)=b_{in}\leq 0$ 
for all $i\in\{1\ld n-1\}$, implying 
$W_\cX\leq\{w\in W';\al_0^{\prime w}=\al_0'\}$.
To show the converse, as above we conclude that 
$\{w\in W';\al_0^{\prime w}=\al_0'\}$
is generated by the reflections it contains,
and any such reflection is of the form $s_{\bt'}\in W'$
for some positive root $\bt'=\sum_{i=1}^{n-1}b_i\al'_i\in\Phi'$,
\ie we have $b_i\geq 0$ for all $i$.
Thus from $B(\bt',\al'_0)=0$, and $B(\al'_i,\al'_0)\leq 0$ for all 
$i\in\{1\ld n-1\}$, we infer $b_i=0$ whenever $B(\al'_i,\al'_0)<0$, 
hence $\bt'\in\lr{\cX}_\R\cap\Phi'$ and thus $s_{\bt'}\in W_\cX$,
implying $W_{\cX}=\{w\in W';\al^{\prime w}=\al'\}$.
\QED

\abs
Thus in this case we have $C_{W'}(s_n)=\{w\in W';\al^{\prime w}=\al'\}$.
Hence the permutation action of $W'$ on the cosets
of $C_{W'}(s_n)$ is isomorphic to its action on 
$(\al')^{W'}\sseq\Phi'$, where $(\al')^{W'}$ encompasses all of $\Phi'$ 
if all roots have the same length, and is the set of all long roots 
if $\Phi'$ contains roots of different lengths.

\AbsC{4. Examples}

\abs
We consider the irreducible finite and affine cases of rank at least $3$,
and a particular irreducible infinite compact hyperbolic case.

\abs
{\bf Type $A_{n-1}$, $n\geq 2$.}
\hfill\includegraphics[width=40mm]{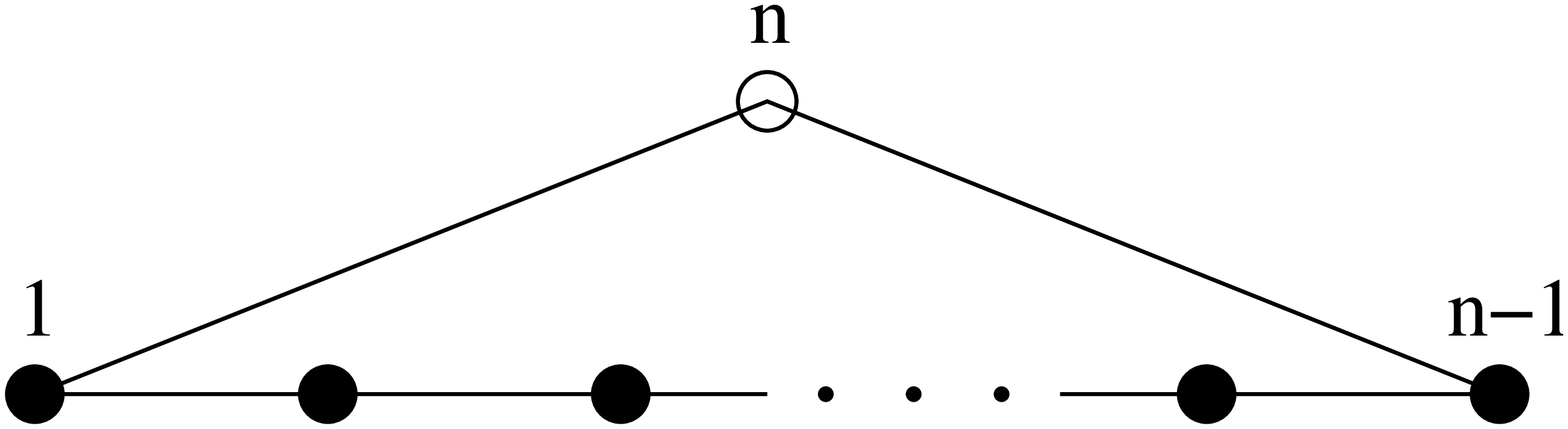} 

\abs
This is the archetypical example of 
symmetric generation \cite[Exc.3.3(1)]{C3}:
We have $W(A_{n-1})\cong\cS_n$, the symmetric group on
$n$ points \cite[Pl.I]{B}, 
where in the natural permutation representation
the distinguished Coxeter generators are the adjacent
transpositions $s_i=(i,i+1)$ for $i\in\{1\ld n-1\}$.

\abs
Letting $t:=s_{n-1}=(n-1,n)$ we get 
$W'=\lr{s_1\ld s_{n-2}}\cong W(A_{n-2})\cong\cS_{n-1}$ and 
$C_{W'}(t)=\lr{s_1\ld s_{n-3}}\cong W(A_{n-3})\cong\cS_{n-2}$, 
where we let $W(A_0):=\cS_{1}=\{1\}$ and $W(A_{-1}):=\cS_0=\{1\}$.
Hence $\cS_n$ is symmetrically generated by 
$(t)^{\cS_{n-1}}=\{(i,n);i\in\{1\ld n-1\}\}$,
where $\cS_{n-1}$ acts on $(t)^{\cS_{n-1}}$
by the natural action on $n$ points:
$\cS_n\cong(2^{\ast(n-1)}\cn\cS_{n-1})/\lr{\!\lr{(s_{n-2}t)^3}\!}$.

\abs
More generally, letting $t:=s_{n-k}=(n-k,n-k+1)$ for 
$k\in\{1\ld\fl{\frac{n}{2}}\}$, we get
$W'=\lr{s_1\ld s_{n-k-1},s_{n-k+1}\ld s_{n-1}}
 \cong W(A_{n-k-1})\tm W(A_{k-1})\cong\cS_{n-k}\tm\cS_k$ and
$C_{W'}(t)=\lr{s_1\ld s_{n-k-2},s_{n-k+2}\ld s_{n-1}}
 \cong W(A_{n-k-2})\tm W(A_{k-2})\cong\cS_{n-k-1}\tm\cS_{k-1}$.
Thus $\cS_n$ is symmetrically generated by 
$(t)^{\cS_{n-k}\tm\cS_k}=\{(i,j);i\in\{1\ld n-k\},j\in\{n-k+1\ld n\}\}$,
where $\cS_{n-k}\tm\cS_k$ acts on $(t)^{\cS_{n-k}\tm\cS_k}$
by the natural action on $(n-k)\cdot k$ points:
$$ \cS_n\cong(2^{\ast (n-k)k}\cn(\cS_{n-k}\tm\cS_k))/\lr{\!\lr{
   (s_{n-k-1}t)^3,\,\,\, (s_{n-k+1}t)^3 }\!} $$

\abs
We summarise this in the following table:
$$ \begin{array}{|l|l|l|l|l|}
\hline
t & W' & C_{W'}(t) & m & \text{relations} \\
\hline
\hline
s_{n-1} & A_{n-2} & A_{n-3} & n-1 & (s_{n-2}t)^3 \\
\hline
s_{n-k} & A_{n-k-1}\tm A_{k-1} & A_{n-k-2}\tm A_{k-2} & k(n-k) & 
                                 (s_{n-k-1}t)^3,(s_{n-k+1}t)^3 \\
\hline
\end{array} $$

\abs
{\bf Type $\wti{A}_{n-1}$, $n\geq 3$.}
Letting $t:=s_n$, we get that 
$W(\wti{A}_{n-1})$ is symmetrically generated by $(t)^{\cS_n}$,
where $W'=W(A_{n-1})\cong\cS_n$ acts on $(t)^{\cS_n}$
by the natural action on ordered pairs
of distinct points in $\{1\ld n\}$,
which coincides with its action on the
roots in the root system of type $A_{n-1}$:
$$ \begin{array}{|l|l|l|l|l|}
\hline
t & W' & C_{W'}(t) & m & \text{relations} \\
\hline
\hline
s_n & A_{n-1} & A_{n-3} & n(n-1) & (s_1 t)^3,(s_{n-1}t)^3 \\
\hline
\end{array} $$

\abs
\includegraphics[width=40mm]{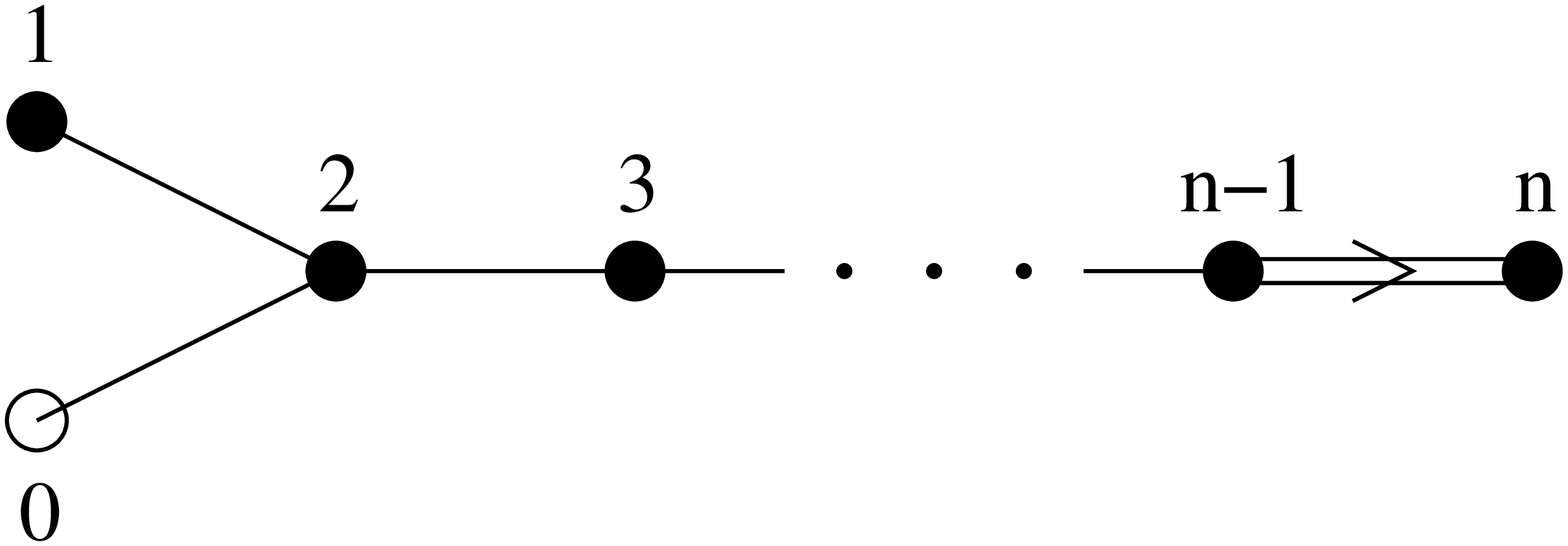}

\vspace*{-3em}
\hfill\includegraphics[width=50mm]{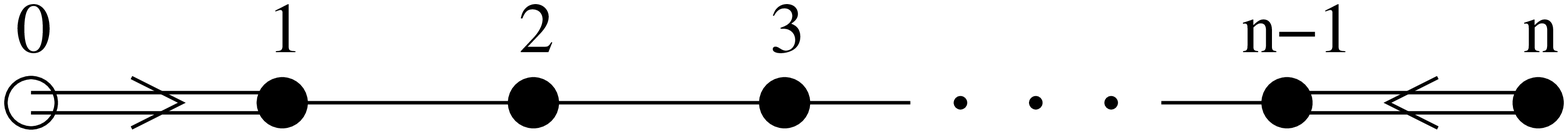}
\vspace*{2em}

\abs
{\bf Types $B_n$, $n\geq 3$, and $C_n$, $n\geq 2$.}
We have $W(B_n)\cong W(C_n)\cong 2^n\cn\cS_n$ \cite[Pl.II, III]{B},
where the geometric representation restricts to 
the natural permutation representation of 
$\lr{s_1\ld s_{n-1}}\cong\cS_n$,
while $s_n$ with respect to the associated permutation basis
acts by the diagonal matrix $\diag[1\ld 1,-1]$.
Hence 
$W(B_n)$ is isomorphic to the group of all monomial 
$\{\pm 1\}$-matrices of size $n$.

\abs
We get the following presentations,
where $k\in\{2\ld n-2\}$, but we only have weak symmetric generation.
To describe the various actions of $W'$ on $(t)^{W'}$
it suffices to note that $W(B_n)\cong W(C_n)$ acts 
on the cosets of $W(B_{n-1})\cong W(C_{n-1})$ as it acts
on the short roots in the root system of type $B_n$,
or equivalently on the long roots in the root system of type $C_n$:
$$ \begin{array}{|l|l|l|l|l|}
\hline
t & W' & C_{W'}(t) & m & \text{relations} \\
\hline
\hline
s_n & A_{n-1} & A_{n-2} & n & (s_{n-1}t)^4 \\
\hline
s_{n-1} & A_{n-2}\tm A_1 & A_{n-3} & 2(n-1) & (s_{n-2}t)^3,(s_n t)^4\\
\hline
s_{n-k} & A_{n-k-1}\tm B_{k} &
          A_{n-k-2}\tm B_{k-1} & 2k(n-k) & (s_{n-k-1}t)^3,(s_{n-k+1}t)^3\\
\hline
s_1 & B_{n-1} & B_{n-2} & 2(n-1) & (s_2 t)^3 \\
\hline
\end{array} $$

\abs
{\bf Type $D_n$, $n\geq 4$.}
\hfill\includegraphics[width=40mm]{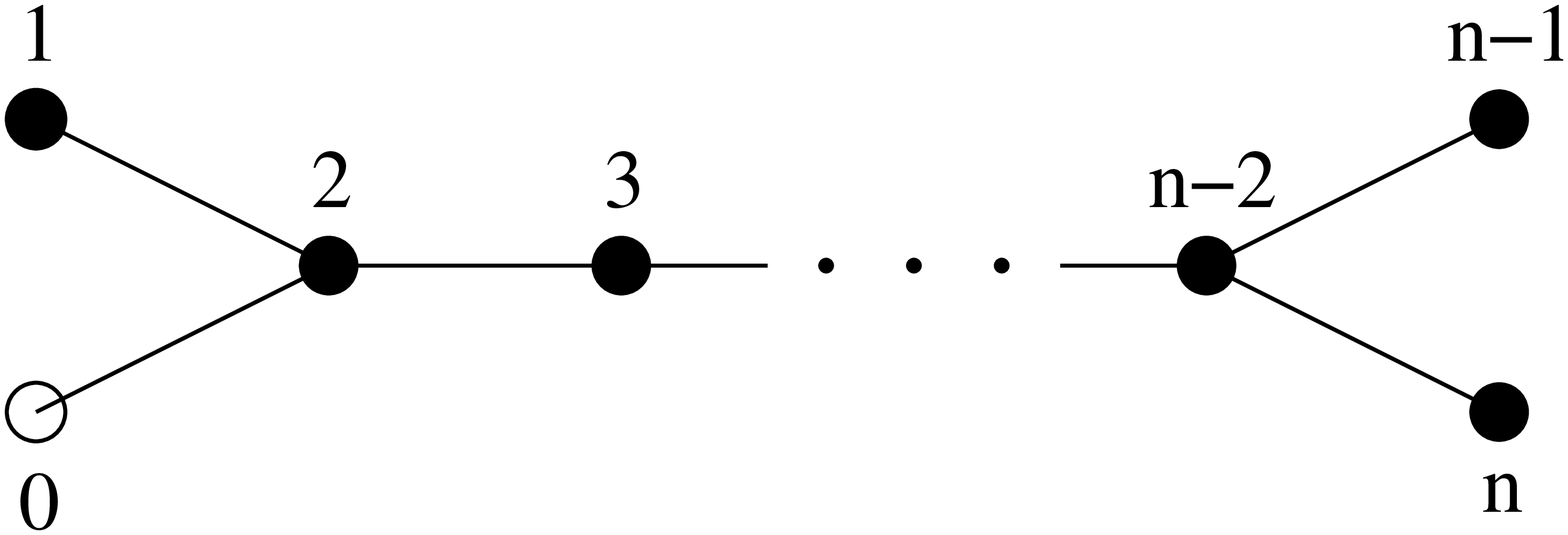} 

\abs
We have $W(D_n)\cong 2^{n-1}\cn\cS_n$ \cite[Pl.IV]{B},
where again the geometric representation restricts to the
the natural permutation representation of 
$\lr{s_1\ld s_{n-1}}\cong W(A_{n-1})\cong\cS_n$.
Moreover, viewing $W(B_n)$ as the group of all monomial 
$\{\pm 1\}$-matrices of size $n$, we have an embedding $W(D_n)<W(B_n)$
such that $W(D_n)$ consists of all such matrices
having an even number of entries $-1$.

\abs
We get symmetric generation subject to the following presentations,
where $k\in\{3\ld n-2\}$ and we let $D_3:=A_3$ and $D_2:=A_1\tm A_1$.
To describe the various actions of $W'$ on $(t)^{W'}$ 
it suffices to note that $\cS_n$ acts on the cosets of $\cS_{n-2}\tm\cS_2$
by the natural action on unordered pairs
of distinct points in $\{1\ld n\}$, and that $W(D_n)$ acts on the cosets 
of $W(D_{n-1})$ as on the short roots in the root system of type $B_n$:
$$ \begin{array}{|l|l|l|l|l|}
\hline
t & W' & C_{W'}(t) & m & \text{relations} \\
\hline
\hline
s_n & A_{n-1} & A_{n-3}\tm A_1 & \binom{n}{2} & (s_{n-2}t)^3 \\
\hline
s_{n-2} & A_{n-3}\tm A_1\tm A_1 & A_{n-4} & 4(n-2) &  
          (s_{n-3}t)^3,(s_{n-1}t)^3, \\
        & & & & (s_n t)^3 \\
\hline
s_{n-k} & A_{n-k-1}\tm D_k & A_{n-k-2}\tm D_{k-1} & 2k(n-k) & 
          (s_{n-k-1}t)^3, \\
        & & & & (s_{n-k+1}t)^3 \\
\hline
s_1 & D_{n-1} & D_{n-2} & 2(n-1) & (s_2 t)^3 \\
\hline
\end{array} $$

\abs
{\bf Type $E_n$, $n\in\{6,7,8\}$.}
\hfill\includegraphics[width=40mm]{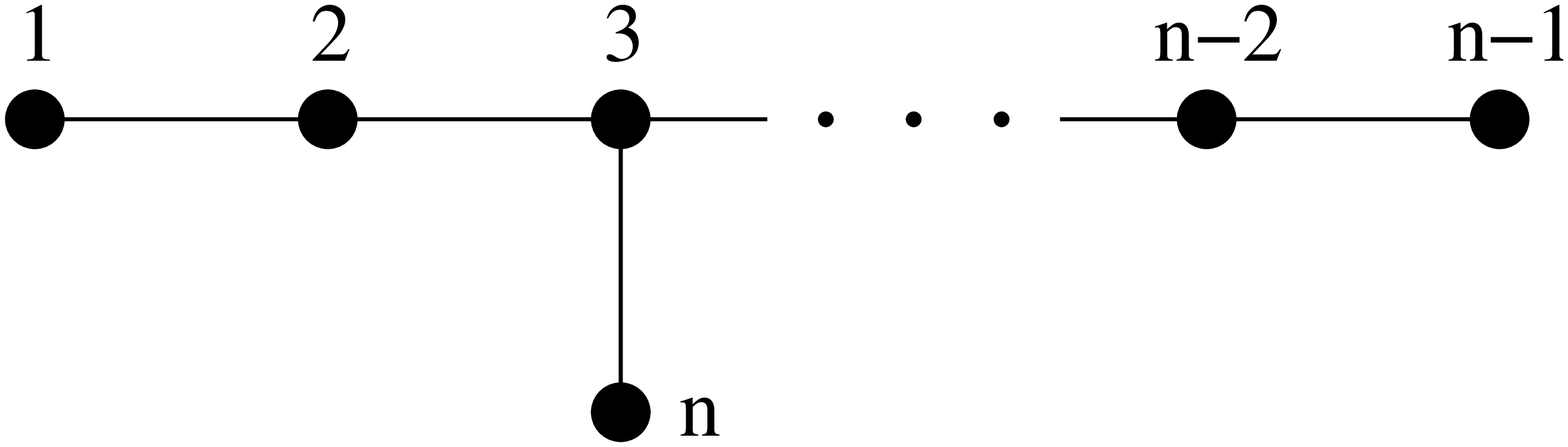}

\abs
We have $W(E_6)\cong GO_6^-(2)=O_6^-(2).2$, 
$W(E_7)\cong 2\tm GO_7(2)=2\tm O_7(2)$ and 
$W(E_8)\cong 2.GO_8^+(2)=2.O_8^+(2).2$, 
\cite[Pl.V--VII]{B} and \cite[p.26, 46, 85]{Atlas}.
We get symmetric generation subject to the following presentations,
where we still let $W(A_0)=W(A_{-1})=\{1\}$.

\abs
To describe the various actions of $W'$ on $(t)^{W'}$ we note 
that $\cS_n$ acts on the cosets of $\cS_2\tm\cS_{n-2}$ and of 
$\cS_3\tm\cS_{n-3}$ by the natural action on unordered pairs and triples 
of pairwise distinct points in $\{1\ld n\}$, respectively,
and that $W(D_n)\cong 2^{n-1}\cn\cS_n$ acts on the cosets of 
$W(A_{n-1})\cong\cS_n$ by its affine action on 
the regular normal subgroup $2^{n-1}$, \ie $\cS_n$ and $2^{n-1}$ 
act by conjugation and translation, respectively.
Moreover, $W(E_6)\cong GO_6^-(2)$ acts on the cosets of $W(D_5)$
as on the isotropic vectors in the associated quadratic space, and
$W(E_7)$ acts on the cosets of $W(E_6)$
as on a certain set of $56$ vectors in the root system of type $E_8$:
$$ \begin{array}{|l|l|l|l|l|}
\hline
t & W' & C_{W'}(t) & m & \text{relations} \\
\hline
\hline
s_n & A_{n-1} & A_2\tm A_{n-4} & \binom{n}{3} & (s_3 t)^3 \\
\hline
s_7, n=8 & E_7 & E_6 & 56 & (s_6 t)^3 \\ 
\hline
s_6, n=7 & E_6 & D_5 & 27 & (s_5 t)^3 \\ 
s_6, n=8 & E_6\tm A_1 & D_5 & 54 & (s_5 t)^3,(s_7 t)^3 \\ 
\hline
s_5, n=6     & D_5 & A_4 & 16 & (s_4 t)^3 \\ 
s_5, n\geq 7 & D_5 \tm A_{n-6} & A_4\tm A_{n-7} & 16(n-5) & 
               (s_4 t)^3, (s_6 t)^3 \\ 
\hline
s_4 & A_4\tm A_{n-5} & A_2\tm A_1\tm A_{n-6} & 10(n-4) & 
      (s_3 t)^3,(s_5 t)^3  \\ 
\hline
s_3 & A_2\tm A_{n-4}\tm A_1 & A_1\tm A_{n-5} & 6(n-3) & 
      (s_2 t)^3,(s_4 t)^3, \\
    & & & & (s_n t)^3 \\
\hline
s_2 & A_1\tm A_{n-2} & A_1\tm A_{n-4} & 2\binom{n-1}{2} & 
      (s_1 t)^3,(s_3 t)^3 \\
\hline
s_1 & D_{n-1} & A_{n-2} & 2^{n-2} & (s_2 t)^3 \\
\hline
\end{array} $$

\abs
{\bf Type $F_4$.} \hfill\includegraphics[width=30mm]{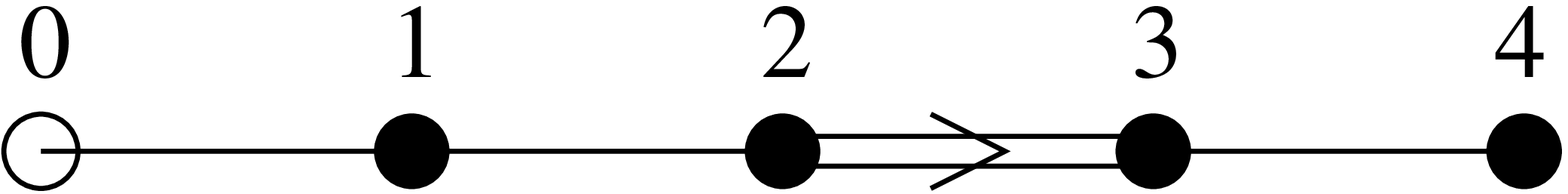} 

\abs
We have $W(F_4)\cong W(D_4)\cn\cS_3\cong (2^3\cn\cS_4)\cn\cS_3$ 
\cite[Pl.VIII]{B}.
We get the following presentations, 
but only weak symmetric generation.
To describe the various actions of $W'$ on $(t)^{W'}$ it suffices to note 
that $W(B_3)\cong 2^3\cn\cS_3$ acts on the cosets of $W(A_2)\cong\cS_3$ 
by its affine action on the regular normal subgroup $2^3$:
$$ \begin{array}{|l|l|l|l|l|}
\hline
t & W' & C_{W'}(t) & m & \text{relations} \\
\hline
\hline
s_4 & B_3 & A_2 & 8 & (s_3 t)^3 \\
\hline
s_3 & A_2\tm A_1 & A_1 & 6 & (s_2 t)^4,(s_4 t)^3 \\
\hline
\end{array} $$

\abs
\mbox{}\hfill\includegraphics[width=30mm]{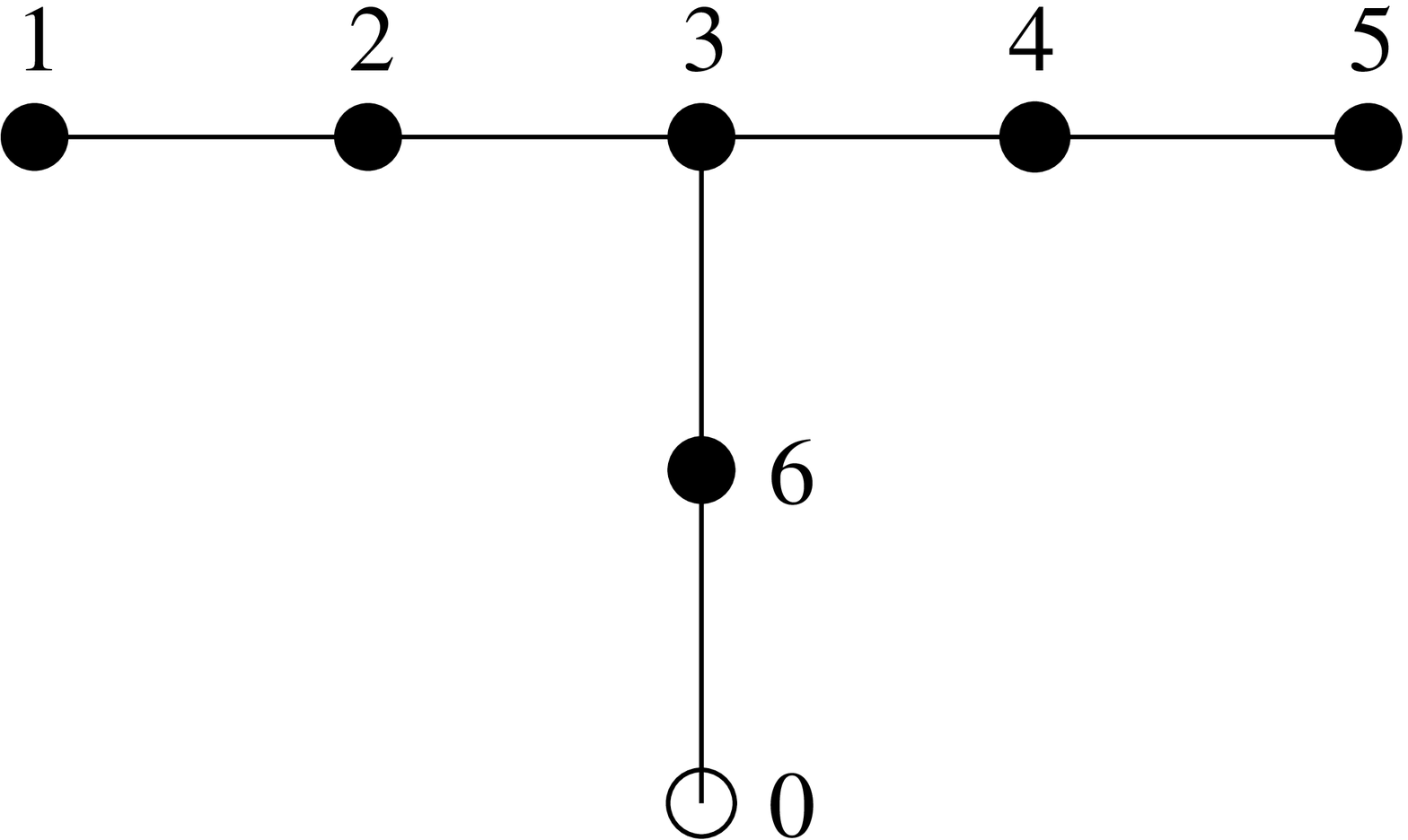}\hfill\mbox{}
\vspace*{-4em}

\abs
\includegraphics[width=45mm]{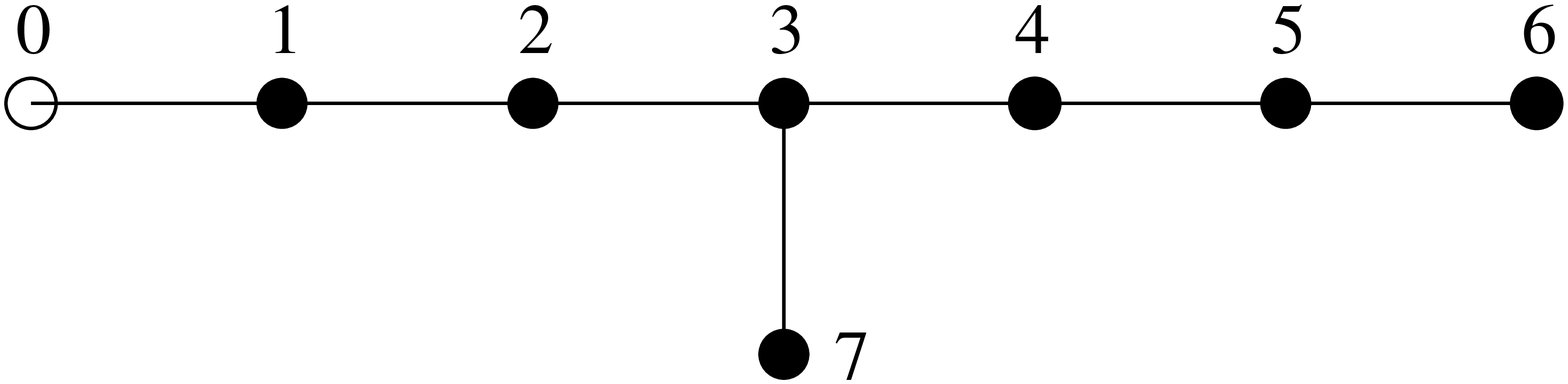}
\hfill\includegraphics[width=50mm]{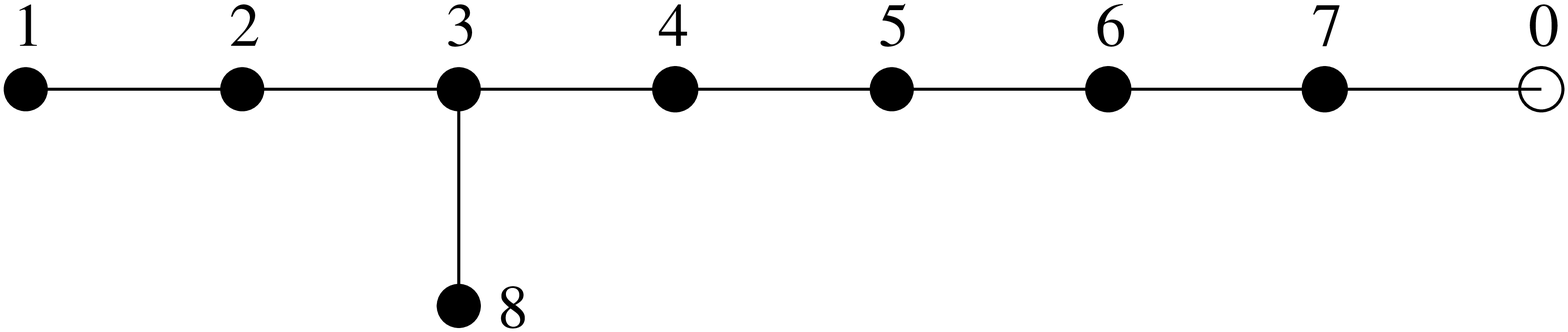}

\abs
{\bf Types $\wti{B}_n$, $\wti{C}_n$, $\wti{D}_n$, $\wti{E}_n$,
$\wti{F}_4$ and $\wti{G}_2$.}
\hfill\includegraphics[width=20mm]{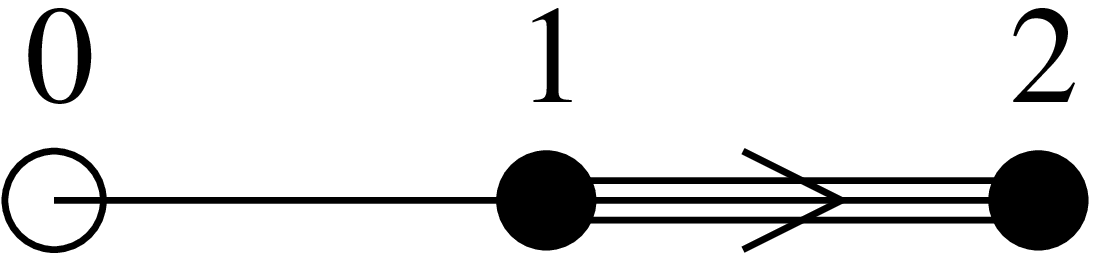}

\abs
Letting $t:=s_0$, we get symmetric generation if and only if the
Dynkin graph associated with $W$ is simply-laced, and 
$W'$ acts on $(t)^{W'}$ as on the long roots in its associated root system:
$$ \begin{array}{|l|l|l|l|l|}
\hline
t & W' & C_{W'}(t) & m & \text{relations} \\
\hline
\hline
s_0 & B_n & A_1\tm B_{n-2} & 2n(n-1) & (s_2 t)^3 \\
s_0 & C_n & C_{n-1} & 2n & (s_1 t)^4 \\
s_0 & D_n & A_1\tm D_{n-2} & 2n(n-1) & (s_2 t)^3 \\
s_0 & E_6 & A_5 & 72 & (s_6 t)^3 \\
s_0 & E_7 & D_6 & 126 & (s_1 t)^3 \\
s_0 & E_8 & E_7 & 240 & (s_2 t)^3 \\
s_0 & F_4 & C_3 & 24 & (s_1 t)^3 \\
s_0 & G_2 & A_1 & 6 & (s_1 t)^3 \\
\hline
\end{array} $$

\abs
{\bf Type $H_n$, $n\in\{3,4,5\}$.}
\hfill\includegraphics[width=35mm]{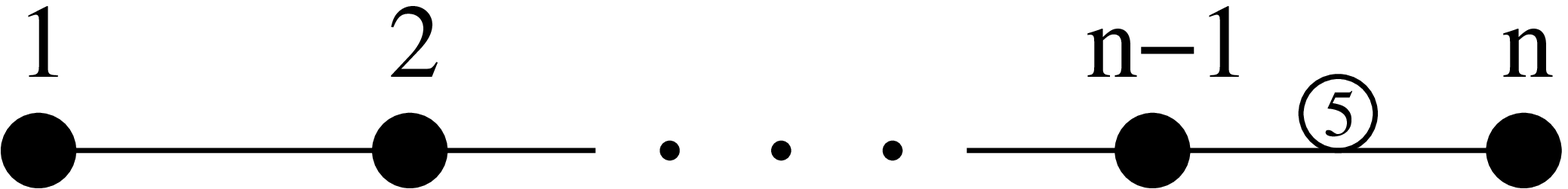}

\abs
We have $W(H_3)\cong 2\tm\cA_5$ and $W(H_4)\cong 2.(\cA_5\tm\cA_5).2$ 
\cite[Exc.VI.4.11, VI.4.12]{B},
where $\cA_5$ denotes the alternating group on $5$ points.
Recall that the geometric representation of the parabolic subgroup
$W(I_2(5))\cong\cD_{10}$ is the group of Euclidean symmetries of 
the regular pentagon. The geometric representation of $W(H_3)$ is 
the group of Euclidean symmetries of the regular dodecahedron,
and hence also of its reciprocal, the regular icosahedron; 
recall that the dodecahedron has pentagonal faces.
The geometric representation of $W(H_4)$ is the group of Euclidean
symmetries of a $4$-dimensional regular polytope having $120$ dodecahedral 
faces \cite[Ch.7.8]{Cox2}.
Finally $W(H_5)$ is an infinite compact hyperbolic group, 
\cite[Exc.V.4.15]{B} and \cite[Ch.6.9]{H}.

\abs
We get symmetric generation subject to the following presentations.
To describe the various actions of $W'$ on $(t)^{W'}$ we note that 
$W(I_2(5))$ acts on the cosets of $W(A_1)$
by the action on the $5$ faces of the regular pentagon,
that $W(H_3)$ acts on the cosets of $W(I_2(5))$
by the action on the $12$ faces of the regular dodecahedron,
and that $W(H_4)$ acts on the cosets of $W(H_3)$  
by the action on the $120$ faces of
the $4$-dimensional regular polytope mentioned above:
$$ \begin{array}{|l|l|l|l|l|}
\hline
t & W' & C_{W'}(t) & m & \text{relations} \\
\hline
\hline
s_n & A_{n-1} & A_{n-2} & n & (s_{n-1}t)^5 \\
\hline
s_{n-1} & A_{n-2}\tm A_1 & A_{n-3} & 2(n-1) & (s_{n-2}t)^3,(s_n t)^5 \\
\hline
s_3, n=5 & A_2\tm I_2(5) & A_1\tm A_1 & 15 & (s_2 t)^3, (s_4 t)^3 \\
\hline
s_2, n=4 & A_1\tm I_2(5) & A_1 & 10 & (s_1 t)^3, (s_3 t)^3 \\
s_2, n=5 & A_1\tm H_3 & I_2(5) & 24 & (s_1 t)^3, (s_3 t)^3 \\
\hline
s_1, n=3 & I_2(5) & A_1 & 5 & (s_2 t)^3 \\
s_1, n=4 & H_3 & I_2(5) & 12 & (s_2 t)^3 \\
s_1, n=5 & H_4 & H_3 & 120 & (s_2 t)^3 \\
\hline
\end{array} $$


\abs
{\sc
B. F.:
School of Mathematics, The University of Birmingham \\
The Watson Building, Birmingham, B15 2TT, United Kingdom} \\
{\sf fairbaib@for.mat.bham.ac.uk}

\abs
{\sc J. M.:
Lehrstuhl D f\"ur Mathematik, RWTH Aachen \\
Templergraben 64, 52062 Aachen, Germany} \\
{\sf Juergen.Mueller@math.rwth-aachen.de}

\end{document}